\documentclass{ieeetran}

\IEEEoverridecommandlockouts                              % This command is only
\usepackage{amsmath} % assumes amsmath package installed
\usepackage{amssymb}  % assumes amsmath package installed
\usepackage{amsfonts}
\usepackage{graphicx}

\newtheorem{thm}{Theorem}

\newtheorem{lem}{Lemma}
\newtheorem{rem}{Remark}

\newtheorem{defn}{Definition}
\newtheorem{ass}{Assumption}

\newenvironment{prf}{\smallbreak\noindent{\it Proof: }}{\hfill$\Box$\smallbreak}

\newcommand{\cl}[1]{{\cal #1}}

\title{\LARGE \bf
Identification of Nonlinear Systems with \\Stable Limit Cycles via Convex Optimization
}

\IEEEoverridecommandlockouts
\author{
Ian R. Manchester \  Mark M. Tobenkin \  Jennifer Wang
\thanks{Ian R. Manchester is with the Australian Centre for Field Robotics, School of Aerospace, Mechanical and Mechatronic Engineering, University of Sydney, NSW, 2006, Australia. Email:ian.manchester@sydney.edu.au}
\thanks{Mark M. Tobenkin is with the Department of Electrical Engineering and Computer Science, Massachusetts Institute of Technology, Cambridge, MA 02139, USA.}
\thanks{Jennifer Wang is with the Center for Human Genetics Research, Massachusetts General Hospital, Boston, MA 02114, USA.}
\thanks{This was supported in part by National Science Foundation Grant No. 0835947.}}

\begin{document}

\maketitle
\begin{abstract}
We propose a convex optimization procedure for black-box identification of
nonlinear state-space models for systems that exhibit stable limit cycles (unforced periodic solutions). It extends the ``robust identification error'' framework in which a convex
upper bound on simulation error is optimized to fit rational
polynomial models with a strong stability guarantee. In this work, we
relax the stability constraint using the concepts of transverse dynamics and orbital stability, thus
allowing systems with autonomous oscillations to be identified. The resulting optimization problem is convex, and can be formulated as a semidefinite program. A simulation-error bound is proved without assuming that the true system is in the model class, or that the number of measurements goes to infinity. Conditions which guarantee existence of a unique limit cycle of the model are proved and related to the model class that we search over. The
method is illustrated by identifying a high-fidelity model from experimental recordings of a live rat hippocampal neuron in culture.
\end{abstract}

\section{Introduction}

Black-box identification of highly nonlinear systems poses many challenges, including flexibility of representation, efficient optimization of parameters, model stability,  and accurate long-term simulation fits \cite{Ljung10, Sjoberg95}. It is especially challenging when the system exhibits autonomous oscillations: such a system is intrinsically nonlinear and lives on the ``edge of stability'', since periodic solutions must have at least one critically-stable Lyapunov exponent \cite{Hale80}. 

Recently, a new framework has been introduced for identifying a broad class of nonlinear systems along with certificates of model stability and accuracy of long-term predictions \cite{Tobenkin10a}. However, this method necessarily fails if the system has autonomous oscillations. In this paper we extend the method of \cite{Tobenkin10a} to remove this restriction.

The main contribution of this paper is a method to identify highly nonlinear systems which:
\begin{itemize}
\item searches over a very broad class of models, including those with limit cycles,
\item guarantees a (local) bound on deviation of open-loop model simulation from the data records,
\item is posed as a convex optimization problem,
\item is analysed without assuming that the true system is in the model class or that the number of measurements grows to infinity.
\end{itemize}

\subsection{Identification of Oscillating Systems}

In many scientific fields there is a
need to capture oscillatory behaviour in the form of a compact
mathematical model which can then be used for simulation, analysis, or control design. When the data comes from experimental recordings, this is known as {\em system identification}. It is also becoming more frequent to perform {\em model-order reduction} via system identification methods from solutions of a very high dimensional simulation, e.g. computational fluid dynamics \cite{Lucia04} or a detailed electronic circuit model (see, e.g., \cite{Bond10, Demir00}).

In biology, systems that oscillate seem to be the rule rather than the exception: heartbeats, firefly synchronization, circadian rhythms, neuron spking, and many others \cite{Rapp87, Murray02}. Nonlinear oscillator models have been used in speech analysis and synthesis, where stability of the identified model has been acknowledged as a major issue \cite{Kubin05}. 

Reduced-order modelling of oscillations in computational fluid dynamics has recently been approached via proper orthogonal decomposition (POD) and Galerkin methods \cite{rowley2000pod, rowley2004model}. It was noted that, although local stability of models can be guaranteed for equilibria by careful choice of projection operators, the same cannot be said for limit cycle solutions. In fact, it was frequently observed that the reduced model would diverge from the target oscillation \cite{rowley2000pod}.

To the authors' knowledge, there is no generally applicable methods of system identification -- or model-order reduction -- for oscillating systems. One family of approaches popular for aerospace model reduction is harmonic balance (describing function) methods, in which the period of oscillation is assumed known and the model is reduced by considering the problem in the Fourier-series domain
\cite{Beran04}, \cite{Kerschen06}, \cite{Lucia04}. A similar approach has been taken to analyse phase-locked loops and oscillators, in which a local phase-offset system is of primary interest \cite{Demir00}.
 Neither of these approaches extend easily to situations in which the frequency of oscillation is input-dependent.
Other papers assume a known decomposition into a stable linear part and a static nonlinear map, and consider it a problem of closed-loop linear system identification \cite{Casas02}. Applications have included identification of combustion instabilities \cite{Murray98, Dunstan01}. A mixed empirical/physics-based approach has been used to produce low-order models of periodic vortex shedding in fluid flows \cite{Noack03}.

\subsection{Stability of Oscillations}

No linear system can produce an asymptotically stable limit
cycle. Identifying nonlinear models from data is a difficult problem,
primarily due to the complex relationship between system
parameters and long term behaviour of solutions. A recent approach,
which this paper builds upon, works via convex optimization of a {\em robust identification
  error} which imposes an asymptotic stability constraint on the identified model \cite{Tobenkin10a}.

However, if the system has a periodic solution, not driven by a periodic forcing term  then
this approach must fail: the stability constraint is too strong. To see this, suppose a system
\[
\dot x = f(x) \in \mathbb R^n
\]
has a non-trivial $T$-periodic solution $x^\star(t)$, then
$x^\star(t+\tau), \tau \in (0, T)$ is also a solution which will never
converge to $x^\star(t)$.

The natural notion of stability for oscillating systems is {\em orbital stability}. A T-periodic solution $x^\star$ is orbitally stable if nearby initial conditions converge to the solution {\em set} in state space $\mathcal X = \{x(\tau): \tau \in [0, T]\}$ and not necessarily to the particular time solution $x^\star(t)$. This is a weaker condition than standard (Lyapunov) asymptotic stability.

Orbital stability can be studied via the introduction of so-called {\em transverse coordinates}, also referred to as the {\em moving Poincar\'e section}
\cite{Hale80, Leonov06}. The basic idea is to construct a new coordinate system at each point of the solution, decomposing the state into a scalar component tangential to the solution curve, and a component of dimension $n-1$ transversal (often orthogonal) to the solution curve. 

It is known that  periodic solution of a nonlinear differential equation is orbitally stable if and only if the dynamics in the transverse coordinates are stable \cite[Chap. VI]{Hale80}. This framework has previously been used to design stabilizing controllers and analyze regions of attraction for oscillating systems \cite{Hauser94a, Shiriaev08, Manchester10a, Manchester10, Shkolnik10a}. It has also been used to analyze the convergence of prediction-error methods when identifying a linear/static-nonlinearity feedback interconnection that can oscillate \cite{Casas02}. In this paper we extended the robust identification error method of \cite{Tobenkin10a} using a storage function in the transverse coordinates, so as to robustly identify a broad class of nonlinear systems that may (or may not) admit autonomous oscillations.

A preliminary version of this paper was presented in  \cite{manchester2011identification}.

\subsection{Paper outline}
The structure of the paper is as follows: in Section \ref{sec:prob} we set up the mathematical problems statement; in Section \ref{sec:rie} we review the method proposed in \cite{Tobenkin10a} and explain why it is not suitable for oscillating systems; in Section \ref{sec:trie} we outline proposed approach and prove the main theoretical results; in Section \ref{sec:relax} we give a convex (semidefinite) relaxation of the associated optimization problem; in Section \ref{sec:contraction} we present a result guaranteeing existence of limit cycles for models, and relate it to the identification algorithm we propose; in Section \ref{sec:impl} we discuss practical matters of implementation and the utility of the model class; in Section \ref{sec:exp} we present experimental results fitting membrane potential dynamics of a spiking rat hippocampal neuron in culture; Section \ref{sec:conc} has some brief conclusions; in two appendices we provide details of the experimental setup and a technical lemma used in the proof of the main result.

%%% Local Variables: 
%%% mode: latex
%%% TeX-master: "transID_cdc.tex"
%%% End: 

\section{Problem Statement}\label{sec:prob}

Given a data record of states, inputs, and outputs $\{\tilde x(t), \tilde u(t), \tilde y(t)\}, t\in [0, T]$, the general problem is to construct a compact model in the form of a differential equation that, when simulated, faithfully reproduces the data. Here we assume that the data record consists of smooth continuous-time signals on an interval, though in practice it will consist of a finite sequence of data points. To pose the problem exactly we must specify both a model class to search over, and an optimization objective.

\subsection{Model Class}

The model class we will search over consists of continuous-time state-space models with state $x\in
\mathbb R^n$, input $u\in\mathbb R^m$, output $y\in\mathbb R^p$, and
dynamics defined in the following implicit form:
\begin{eqnarray}
\frac{d}{dt}e(x)&=&f(x,u), \label{eqn:model_ef}\\
y &=&g(x,u),\label{eqn:model_g}
\end{eqnarray}
where $e: \mathbb R^n \rightarrow \mathbb R^n, f: \mathbb R^n \times
\mathbb R^m \rightarrow \mathbb R^n, g: \mathbb R^n \times
\mathbb R^m \rightarrow \mathbb R^p$ are smooth functions. 
The Jacobians with respect to $x$ of $e(x)$, $f(x,u)$, and $g(x,u)$ are denoted
$E(x) = \frac{\partial}{\partial x}e(x), \ F(x,u) = \frac{\partial}{\partial x}f(x,u), G(x,u) = \frac{\partial}{\partial x}g(x,u)$.
We will enforce the constraint that $E(x)$ be nonsingular, so the above implicit model can equivalently be written in explicit form:
\[
\dot x=E(x)^{-1}f(x,u).
\]

\begin{rem} To implement the methods described in this paper, $e(x), f(x,u),$ and $g(x,u)$ should come from a finite-dimensional convex class of matrix functions for which one can efficiently check positivity. In practice, we use matrices of polynomials or trigonometric polynomials and make use of the sum-of-squares relaxation to prove positivity \cite{Parrilo03a, Megretski03}.
\end{rem}

\subsection{Optimization Objective}

The general problem we consider is to minimize, over choice of $e, f, g$, the value of the {\em simulation error:}
\[
\frak{E}= \int_0^T|y(t)-\tilde y(t)|^2dt
\]
where $y(t)$ is the solution of \eqref{eqn:model_ef},
\eqref{eqn:model_g} with $x(0) = \tilde x(0)$. One may also wish to ensure that the dynamical system defined by \eqref{eqn:model_ef}, \eqref{eqn:model_g} is well-posed and has some sort of stability property. Note that we do not assume that the system from which data is recorded is in the model class.

Direct optimization of simulation error is not usually tractable: the relationship between system parameters and model simulation is highly nonlinear, and for black-box models we typically don't have good initial parameter guesses. We make the problem tractable (a convex program) through a series of approximations and relaxations.

A further problem arises when the system exhibits a limit cycle; namely, even if the system is modelled very accurately, but the initial condition has a small error in phase, then the simulation error can be very large. A similar problem is caused by inaccuracies in the phase dynamics, which are expected when the true system is not in the model class.

A common and straightforward approach to approximating dynamics is to minimize {\em equation error} by basic least squares, i.e. to minimize
\[
\sum_i |E(\tilde x(t))\dot {\tilde x}(t)-f(\tilde x(t),\tilde u(t))|^2+|\tilde y(t)-g(\tilde x(t), \tilde u(t))|^2
\]
or a similar criterion\footnote{For the implicit models we consider, a further constraint is needed to prevent $E(x)=0$ and $f(x,u)=0$ being the optimum, e.g. a well-posedness condition on $E(x)$, which we discuss later.}. The advantage is that the optimal solution can be computed extremely efficiently by solution of a linear system. The disadvantage is that minimizing equation error gives no guarantees about long-term simulation of the model, nor even that the model is stable. For highly nonlinear systems such as those exhibiting limit cycles, this is especially problematic.

% We will also consider the linearised system about a trajectory:
% \[
% E(x)\dot\Delta = F(x,u)\Delta
% \]
% where $F(x,u)$ is the Jacobian of $f(x,u)$ with respect to $x$. 

% \subsection{Orbital Stability}

% \begin{defn}\label{defn:orbital_stab}
% A model of class \eqref{eqn:model_ef}, \eqref{eqn:model_g} is
% considered orbitally stable if for any two initial conditions $x_1(0)$
% and $x_2(0)$ such that $|x_1(0)-x_2(0)|$ is sufficiently small, there
% exists a continuous map $\tau:[0, \infty)\rightarrow [0,
% \infty)$ such that the solutions satisfy
% \[
% |x_1(t)-x_2(\tau(t))| \rightarrow 0
% \]
% \end{defn}

% Equivalent to classical definition for periodic trajectories, stronger
% for non-periodic (Zhukovsky).

%%% Local Variables: 
%%% mode: latex
%%% TeX-master: "transID_cdc"
%%% End: 

\section{Nonlinear System Identification via Robust Identification Error}\label{sec:rie}

\label{sec:rie}
The papers \cite{Tobenkin10a, manchester2012stable} provided ``local'' and ``global'' bounds for simulation error via discrete and continuous-time models. In this section we briefly recap the local results for continuous-time systems and explain why they cannot be directly applied to model systems with autonomous oscillations. The basic idea is to search jointly for system equations as well as a storage function with output reproduction error as a supply rate. Standard dissipation inequality arguments \cite{Willems72} then provide a bound on long-term simulation error.

\subsection{Linearized Simulation Error}

Suppose we have a model of the form \eqref{eqn:model_ef}, \eqref{eqn:model_g} and a data record $\{\tilde x(t), \tilde u(t), \tilde y(t)\}, t\in [0, T]$. We introduce the {\em linearized simulation error} as a local measure of the model's divergence from the data.

First, we define the {\em equation error} signals associated with  \eqref{eqn:model_ef}, \eqref{eqn:model_g} and the data:
\begin{eqnarray}
\epsilon_x(t) &=& E(\tilde x(t))\dot {\tilde x}(t)-f(\tilde x(t),\tilde u(t)), \\ \epsilon_y(t) &=& \tilde y(t)-g(\tilde x(t), \tilde u(t)).
\end{eqnarray}
Now, consider the following family of systems parametrized by $\theta\in[0,1]$:
\begin{eqnarray}
E(x)\dot x &=& f(x,u)+f_\theta,\\
y &=& g(x,u)+g_\theta.
\end{eqnarray}
Let $(x_\theta, y_\theta)$ be the solution of the above system with $f_\theta = (1-\theta) \epsilon_x$ and $g_\theta = (1-\theta) e_y$. That is, for $\theta =1$ we have $x_\theta=x, y_\theta = y$ and for $\theta=0$ we have $x_\theta=\tilde x, y_\theta=\tilde y$.
We can consider the following linearized simulation error about the recorded trajectory:
\[
\mathcal{E}= \lim_{\theta \rightarrow 0} \frac{1}{\theta^2}\int_0^T|y_\theta(t)-\tilde y(t)|^2dt
\]
as local approximation of the true simulation error $\frak E$.

\subsection{Robust Identification Error}

Note that $\mathcal{E}$ can alternately be represented as
\[
\mathcal{E}=\int_0^T|G(\tilde x(t), \tilde u(t))\Delta+\epsilon_y(t)|^2dt
\]
where
\begin{equation}\label{eqn:Delta_t}
\Delta(t) = \lim_{\theta \rightarrow 0} \frac{1}{\theta}[x_\theta(t)-\tilde x(t)]
\end{equation}
which obeys the dynamics
\[
\frac{d}{dt}(E(\tilde x(t))\Delta(t)) = F(\tilde x(t), \tilde u(t))\Delta(t)+\epsilon_x(t).
\]
That is, 
$\Delta(t)$ is an estimate of the deviation of the model simulation $x(t)$ from the recorded data trajectory $\tilde x(t)$.

It was shown in \cite{Tobenkin10a} that
\begin{equation}\label{eqn:rie_bound}
\mathcal{E}\le \int_0^T \bar{\mathcal {E}}_Q(t)dt
\end{equation}
for any $Q=Q'>0$, where\footnote{Here, and frequently throughout the paper, we drop the arguments on $E(\tilde x(t)), F(\tilde x(t), \tilde u(t)), G(\tilde x(t), \tilde u(t)), \epsilon_x(t)$, and $\epsilon_y(t)$ for the sake of compactness of notation. It should be understood that these are always functions of time and the data.}
\begin{equation}\label{eqn:rie_diss}
\bar{\mathcal {E}}_Q(t) = \sup_{\Delta\in\mathbb R^n}\{2\Delta' E'Q(F\Delta+\epsilon_x)+|G\Delta+\epsilon_y|^2\}.
\end{equation}
The systems theory interpretation of \eqref{eqn:rie_diss} is that the first term in the supremum is the derivative of a positive-definite storage function with respect to linearized simulation error, and the second term is the output reproduction error.

The bound \eqref{eqn:rie_bound} suggests searching over functions $e, f, g$ and a matrix $Q=Q'>0$ so as to minimize the right-hand-side of \eqref{eqn:rie_bound}. This optimization is still non-convex, but a convex relaxation is given in \cite{Tobenkin10a} (we use a similar relaxation in Section \ref{sec:relax} of the present paper).

Each of the supremums over $\Delta$ in \eqref{eqn:rie_diss} are finite if and only if the matrices
\[
R= E'QF+F'QE+G'G
\]
for each data point is negative semidefinite. If this property holds for all $x, u$, then it has been proven that the system is globally incrementally output stable. The reason is that $\Delta'E'QE\Delta$ is a contraction metric for the system \cite{lohmiller1998contraction}, and $\Delta'(E'QF+F'QE)\Delta$ is its derivative. A formal proof of stability is given in \cite{Tobenkin10a}.

For the purposes of the present paper, it is sufficient to note that enforcing global incremental stability is {\em too strong} to allow identification of systems exhibiting autonomous oscillations, since such systems {\em cannot} satisfy this property. The main purpose of this paper is to overcome this limitation via a reformulation of the RIE in the transverse dynamics.

%%% Local Variables: 
%%% mode: latex
%%% TeX-master: "transID_cdc"
%%% End: 

\section{Transverse Robust Identification Error}\label{sec:trie}

In oscillating systems, perturbations in phase cannot be stable and will therefore accumulate over time. The natural form of stability is {\em orbital} stability, which can be defined as stability to a solution set in state space, rather than a particular time solution. A standard framework for anlaysis of orbital stability is via transverse coordinates.

Correspondingly, if both a true system and an identified model admit autonomous oscillations, then it is not possible to ensure that phase deviations between them converge in time.  In this section, we adapt the transverse dynamics approach to the problem of system identification, and construct the {\em Transverse Robust Identification Error} (TRIE).

Let $R_T$ denote the class of time reparametrizations, i.e. smooth monotonically increasing  (and hence diffeomorphic) functions $\tau: [0, T]\rightarrow [0, T_\tau]$ for some $T_\tau=\tau(T)>0$.

Therefore we introduce the concept of {\em orbital simulation error}:
\[
\frak E_\tau = \int_0^T|y(\tau(t))-\tilde y(t)|^2d\tau
\]
defined for a particular time reparametrization $\tau(t): [0, T] \rightarrow [0, T_\tau]$ for some $T_\tau$. Note that $T_\tau$ may be greater or less than $T$, depending on whether the simulation model ``leads'' or ``lags'' the data, in this case the simulation model would be run for a longer or shorter time.

Similarly, we define the {\em orbital linearized simulation error}:
\[
\mathcal{E}_\tau = \lim_{\theta \rightarrow 0} \frac{1}{\theta^2}\int_0^T|y_\theta(\tau(t))-\tilde y(t)|^2dt
\]
as local approximation of $\frak E_\tau$.

Note that $\mathcal{E}$ can alternately be represented as
\[
\mathcal{E}=\int_0^T|G(\tilde x(t), \tilde u(t))\Delta(t, \tau(t))+\epsilon_y(t)|^2dt
\]
where
\begin{equation}\label{eqn:Delta_t}
\Delta(t, \tau) = \lim_{\theta \rightarrow 0} \frac{1}{\theta}[x_\theta(\tau)-\tilde x(t)].
\end{equation}

One can consider the ``optimal'' time reparametrization $\tau^\star = \arg\min_{\tau \in R_T} \mathcal E_\tau$. We also make the following assumption on a sub-optimal but computationally tractable alternative:
\begin{defn}
Define the family of transversal surfaces $S(t)$ for $t\in[0, T]$ as
\[
S(t) := \{\Delta: \Delta'\dot{\tilde x}(t) = 0\},
\]
\end{defn}
i.e. the $(n-1)$--dimensional subspace orthogonal to the data vector field at time $t$.

\begin{figure} \centering \includegraphics[width=\columnwidth]{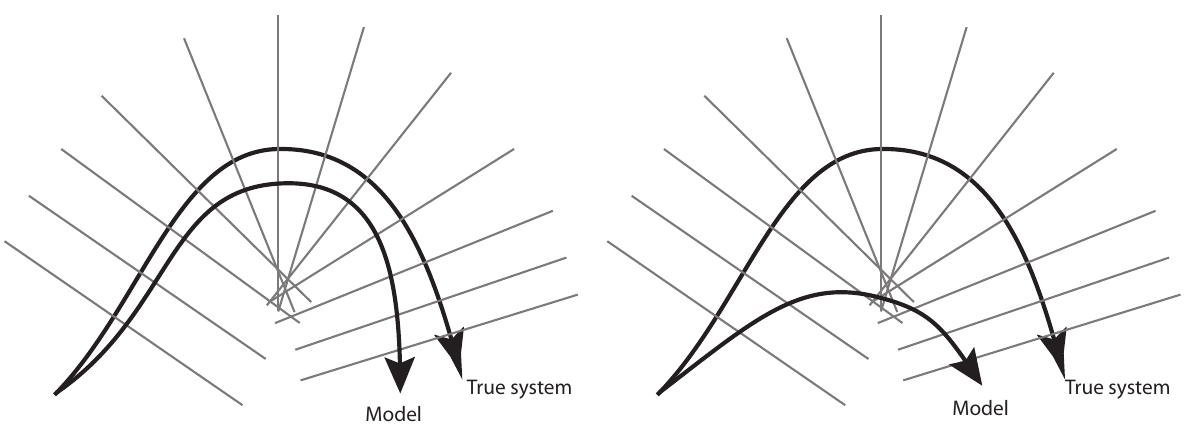}
\caption{Left: an illustration of the solution of a model passing monototonically through the transversal surfaces $S(t)$ orthogonal to the dynamics of the true system $\dot{\tilde x}$, satisfying Assumption \ref{ass:monotonic}. Right: an illustration of a model which does not pass monotonically through the surfaces $S(t)$.}
\label{fig:ts}
\end{figure}

\begin{ass}\label{ass:monotonic}
Assume the linearized simulation error is sufficiently small that there exists a time reparametrization $\tau_S\in R_T$ such that $x(\tau_S(t))\in S(t)$ for all $t\in [0, T]$.
\end{ass}
That is, the linearized simulation error passes monotonically through each of the transversal surfaces $S(t)$ for $t\in[0, T]$. This is illustrated in Figure \ref{fig:ts}.

\begin{ass} For the main theoretical results in this paper, we assume that the input is a constant signal $u(t) = u$ for all $t$.
\end{ass}
Note that the algorithms  we propose can be applied with a time-varying input, but since orbital simulation error explicitly allows the model to be at a different phase of a limit cycle -- and hence different point in state space -- than the true system, it cannot be guaranteed that applying a time-varying input will have the same effect on both the true system and the model, even if the model is perfect.

% Let $\tau_\theta(t)$ be a smooth function such that for $t\in [0,T]$
% \begin{eqnarray}
% \tilde x(\tau_\theta(t))-x_\theta(t) \in S(t).
% \end{eqnarray}
% Observe that $\tau_1(t)=t$ and $\tau_\theta(t)$ is well-defined and smooth function of $\theta$ when $\theta$ is sufficiently close to zero. We will never explicitly compute $\tau_\theta$.

\subsection{Simulation Error Bound}

We are now in a position to give the first main theoretical result of the paper.

\begin{defn}
Define the following projection operators:
\begin{equation}\label{eqn:pi}
\pi(t):= \frac{\dot {\tilde x}(t) \dot  {\tilde x}(t)'}{|\dot  {\tilde x}(t)|^2}, \
\Pi(t):= I-\pi(t),
\end{equation}
i.e. $\pi(t)$ projects on to the one-dimeonsional subspace spanned by $\dot  {\tilde x}(t)$ and $\Pi(t)$ projects on to the $(n-1)$-dimensional subspace orthogonal to this. Furthermore, let $\Pi^r(t)\in\mathbb R^{n\times (n-1)}$  be a matrix with orthonormal columns spanning the subspace orthogonal to $\dot {\tilde x}(t)$, i.e. a ``reduced'' form of the rank $(n-1)$ matrix $\Pi(t)$ containing only independent columns.
\end{defn}

Let 
\begin{eqnarray}
\bar{\cl E}^\perp_Q(t) &=& \sup_{\Delta\in\mathbb R^{n-1}}\{2\Delta{\Pi^r}' E'Q((F+E\dot \Pi) \Pi^r\Delta+\epsilon_x)\notag\\ &&+|G\Pi^r\Delta+\epsilon_y|^2\}\label{eqn:Ebarperp}
\end{eqnarray}
where  and $Q$ a symmetric positive-definite $n\times n$ matrix. 

Note that the above supremums are finite if and only if 
\[
2{\Pi^r}' E'Q((F+E\dot \Pi) \Pi^r+{\Pi^r}'G'G\Pi^r\le 0.
\]
The supremums can be made ``robust'' by enforcing strict negative-definiteness in the above inequality.

\ 

\begin{thm}\label{thm:trie_bound}
Suppose Assumptions 1 and 2 hold. Consider measurement $\tilde x(t), \tilde y(t)$ and simulation $x(t), y(t)$ with the same initial condition $\tilde x(0)=x(0)$ and the same input $u(t)$. Then there exists a time reparametrization $\tau(t)$ such that the following relation holds:
\begin{equation}\label{eqn:E_upperbnd}
\int_0^{T_\tau}|G\Delta(t,\tau)+\epsilon_y|^2d\tau \le \int_0^{T_\tau} \bar{\cl E}^\perp_Q \dot\tau dt.
\end{equation}
\end{thm}

\begin{prf}
The inequality is shown via a dissipation inequality for the following storage function:
\begin{equation}\label{eqn:trans_storage}
 V(\Delta,t) = |E(\tilde x(t)) \Pi(t)\Delta|_Q^2+|\pi(t)\Delta|^2.
\end{equation}

The particular (suboptimal) choice of time reparametrization $\tau(t)$ we consider is that which keeps $\Delta(t,\tau)$ in the surface $S(t)$, i.e. $\Delta(t,\tau)'\dot{\tilde x}(t)=0$.
This choice has two useful properties.

1) $\Pi(t)\Delta(t,\tau) = \Delta(t,\tau)$ and $\pi(t)\Delta(t,\tau) = 0$, by construction of $\Pi(t)$ and $\pi(t)$.

2) For a fixed $t$, the chosen $\tau(t)$ has the property that it minimizes $\frac{d}{d\tau} V(\Delta(t,\tau),t)$ over choices of $\dot \tau$, since the curve $\tilde x(t)$ is orthogonal to $\Delta$.

Using Lemma \ref{lem:Vdot} (see appendix), we see that with $\dot\tau=1$
\begin{equation}\label{eqn:vDbar_dot1}
\frac{d}{d\tau}V(\Delta(t),t)+|G\Delta(t)+\epsilon_y|^2\le  \bar{\cl E}^\perp_Q(t) 
\end{equation}
Due to fact 2) above, it follows that this holds also with the chosen $\dot\tau$ defined by $S_T$.

Since $\tilde x(0)=x(0)$ we have $V(\Delta(0),0)=0$, so integrating \eqref{eqn:vDbar_dot1} gives
\begin{equation}
V(\Delta(T),T)+\int_0^{T\tau}|G\Delta(t)+\epsilon_y|^2d\tau\le \int_0^{T_\tau} \bar{\cl E}^\perp_Q(t)  d\tau.\notag
\end{equation}
and by definition $V(\Delta(T,T_\tau),T)\ge 0$, so the above inequality implies \eqref{eqn:E_upperbnd}. This completes the proof of the theorem.

\end{prf}

\begin{rem} From a data set alone one can integrate $\cl E(t)$ with respect to $\tau$ but not $t$, without computing the model solutions. So in fact we are optimizing with an unknown positive ``weighting'', i.e.
 \begin{equation}
\int_0^T|G\Delta+\epsilon_y|^2dt\le \int_0^{T}\bar{\cl E}^\perp_Q \dot\tau dt.
\end{equation}
Note that if the model is close to the true system, then $\dot\tau\approx 1$, see \cite{Leonov06, Manchester10}, so the weighting factor will not have a great effect.
\end{rem}

%%% Local Variables: 
%%% mode: latex
%%% TeX-master: "transID_TAC"
%%% End: 

\section{A Convex Upper Bound}\label{sec:relax}

Theorem \ref{thm:trie_bound} suggests minimizing the
\[
\int_0^T \cl E^\perp_Q(t) dt
\]
over choices of $e, f, g$ and $Q$ as an effective procedure for system identification. However, this is still a nonconvex optimization. In this
section we propose a convex upper bound for which one can efficiently
find the global minimum via semidefinite programming.

The basic idea is to decompose the each non-convex term into the sum
of a convex and a concave part, and upper-bound the concave part with
a linear relaxation.

\begin{thm}\label{thm:relax}
Define the following quantities,  each of which is linear in the decision variables $e, f$, and $g$.
\begin{eqnarray}
\Delta_e^+&=&E(\tilde x)(I+\dot \Pi)\Pi^r\Delta+F(\tilde x,\tilde u)\Pi^r\Delta+\epsilon_x\notag\\
\Delta_e^-&=&E(\tilde x) (I-\dot \Pi)\Pi^r\Delta-F(\tilde
x,\tilde u)\Pi^r\Delta-\epsilon_x\notag\\
\Delta_y &=&G\Pi^r\Delta+\epsilon_y.\notag
\end{eqnarray}
Then $\bar{\cl E}^\perp_Q(t) \le \hat{\cl E}^\perp_Q(t)$ where
\[
 \hat{\cl E}^\perp_Q
=\sup_{\Delta\in\mathbb R^{n-1}}\left\{\frac{|\Delta_e^+|_Q^2+|\Pi^r\Delta|_{Q^{-1}}^2}{2}-(\Pi^r\Delta)'\Delta_e^-
  + |\Delta_y|^2\right\}
\]
which is convex in $e, f, g,$ and $Q^{-1}>0$.
\end{thm}
\begin{prf}
A similar statement was proved in \cite[Section
V]{Tobenkin10a}. The upper bound is based on the expansion
\[
|a-Q^{-1}\Delta|_Q^2=\Delta Q^{-1}\Delta-2\Delta'a+a'Qa
\]
which, when $Q>0$, clearly implies
\begin{equation}\label{eqn:relaxation}
-a'Qa\le \Delta'Q^{-1}\Delta -2\Delta'a.
\end{equation}
Notice that the right-hand side of \eqref{eqn:relaxation} is convex in
$a$ and $Q^{-1}$ whereas the left-hand-side is concave.

Note that we also have the following expansion:
\begin{equation}\label{eqn:relax2}
4(E\Pi^r\Delta)'Q[(F+E\dot\Pi)\Pi^r\Delta+\epsilon_x]=|\Delta_e^+|_Q^2-|\Delta_e^-|_Q^2.
\end{equation}
The first term on the right-hand-side of \eqref{eqn:relax2} is convex in $e, f$, and
$Q^{-1}$ and the second term is concave. Setting $a=\Delta_e^-$ and applying \eqref{eqn:relaxation} to
\eqref{eqn:relax2} gives the statement of the theorem.
\end{prf}

\subsection{TRIE as an upper bound for equation error}

The results presented so far control the divergence of the model from the data in a ``transversal'' direction, but not in the ``tangential'' or ``phase'' direction. I.e., we have not proven that $\dot\tau \approx 1$. The phase dynamics of a periodic solution acts like a pure integrator \cite{guckenheimer1975isochrons}, so to control simulation error one must simply control equation error in the direction parallel to $\dot {\tilde x}(t)$.

Another reason to want to control equation error directly is that a premise of the above theorems is that the model is already ``good enough'' that it passes monotonically through the transversal surfaces $S(t)$. 

Taking $\Delta=0$, we have
\[
\hat{\cl E}^\perp_Q=\frac{1}{2}|\epsilon_x|_Q^2+|\epsilon_y|
\]
hence TRIE, by taking the supremum over $\Delta$, does also penalize equation error, weighted by $Q$.

One rather artificial circumstance in which this does not result in a bound on equation error occurs when $\dot{\tilde x}(t)$ has the same direction for all $t$. In this case, $\Pi^r$ in the term $|\Pi^r\Delta|_{Q^{-1}}^2$ in $\hat{\cl E}^\perp_Q$, there is nothing in the constraints that stops $Q$ growing in such a way that $\dot{ \tilde x}'Q\dot{ \tilde x}$ is very large, and hence the weighting in the above equation error being very small. However, in practical cases this will never occur.

\subsection{The Proposed Objective for Optimization}
Summarizing the results of this and the previous sections, we have the relations
\[
\frak E^\perp \approx \mathcal E^\perp \le \int_0^T \dot\tau\bar{\cl E}^\perp_Q(t)dt \approx \bar{\cl E}^\perp_Q(t)dt \le \int_0^T \hat{\cl E}^\perp_Q(t)dt,
\]
with the leftmost term being the true orbital simulation error, and the rightmost term being convex in the system equations and $Q^{-1}$. Therefore we propose to perform system identification via the optimization
\[
\int_0^T \hat{\cl E}^\perp_Q(t)dt \rightarrow \min
\]
over choices of $e, f, g,$ and $Q^{-1}>0$ subject to the constraint
$E(x)' + E(x)>I$ for all $x$.

In practice, we will have a record of the true system at a finite number of times $t_i, i = 1, 2, ..., N$, and as a surrogate for the above we minimize the finite sum of the TRIE terms:
\[
\sum_{i=1}^N \hat{\cl E}^\perp_Q(t_i) \rightarrow \min.
\]
Again, the solutions are made ``robust'' by enforcing strictness of the LMI constraints leading finite supremums of each $\hat{\cl E}^\perp_Q(t_i)$

%%% Local Variables: 
%%% mode: latex
%%% TeX-master: "transID_cdc"
%%% End: 

\section{On the Existence of Stable Limit Cycles for Identified Models}\label{sec:contraction}

In \cite{Tobenkin10a} it was proven that if the RIE condition holds everywhere, then the model is globally incrementally $L^2$ stable. This is too strong for systems with limit cycles, but it would be very useful to be able to guarantee that our model has stable limit cycles. In this section we show that the method proposed above will guarantee this property if the error is sufficiently small.

This result is a generalization of the results of \cite{manchester2012contraction} to implicit systems of the form \eqref{eqn:model_ef}. Motivated by contraction theory, we introduce the dynamics of differentials $\delta_x\in\mathbb R^n$, via a linearization of \eqref{eqn:model_ef}:
\[
\frac{d}{dt}(E(x)\delta_x) = F(x,u)\delta_x.
\]
Note that despite the appearance of linearized dynamics, the statements in this section are rigorous results for solutions the true nonlinear system \eqref{eqn:model_ef}, and not based on a local approximation.

\begin{defn}
For $v\in\mathbb R^n$ define $\Pi^r_v$ to be a $\mathbb R^{n\times (n-1)}$ matrix spanned by columns orthogonal to $v$. E.g. for $v=\dot {\tilde x}(t)$ one can take $\Pi^r_{\dot {\tilde x}(t)}=\Pi^r(t)$ from Section \ref{sec:trie}. Similarly, define $\Pi_{\dot x}\in\mathbb R^{n\times n}$ to be an orthonormal matrix projecting on to the subspace orthogonal to $v$.
\end{defn}

\begin{defn}
A compact set $K\subset \mathbb R^n$ is defined to be {\em strictly forward invariant} with respect a dynamical system if it has non-empty interior and any solution starting with $x(0)$ on the boundary of $K$ has $x(t)$ in the interior of $K$ for all $t>0$.
\end{defn}

\begin{thm}\label{thm:contraction}
Suppose that there exists a path-connected set $K$ which is forward-invariant with respect to $E(x)\dot x = f(x,u)$ such that
\begin{equation}\label{eqn:contraction}
\delta_x{\Pi^r_{\dot x}}' E(x)'Q((F(x,u)+E(x)\dot \Pi_{\dot x}) \Pi^r_{\dot x}\delta_x<0
\end{equation}
for all $x\in K$ and for all $\delta_x\in\mathbb R^{n-1}$, where $\dot x=E(x)^{-1}f(x,u)$ is the model derivative. Then there exists a unique periodic solution of the model $x^\star(\cdot)$ in $K$ and from every initial condition $x(0)\in K$, the solutions of $E(x)\dot x = f(x,u)$ converge to the orbit of $x^\star$.
\end{thm}

\begin{prf}We begin by showing that any two solutions converge under possible time reparametrization, i.e. given any two points $x_a, x_b$ in $K$, there exists monotonic smooth functions $\tau_a(t):[0,\infty)\rightarrow [0, \infty)$ and  $\tau_b(t):[0,\infty)\rightarrow [0, \infty)$ such that $x(\tau_a(t))\rightarrow x(\tau_b(t))$ as $t\rightarrow \infty$.

Let $V(x,\delta_x)$ be a positive-definite function of $x$ and $\delta_x$, quadratic in $\delta_x$. Consider two points $x_a$ and $x_b$ in $K$ and a smooth path $\gamma_0:[0,1] \rightarrow K$ connecting them, i.e. $\gamma(0) = x_a, \gamma(1) = x_b$. Define the following measure of distance along the path:
\[
D_\gamma(x_a,x_b) = \int_0^1 \sqrt{V\left(\gamma_0(s), \frac{d\gamma_0 (s)}{ds}\right)}ds.
\]
Clearly if $\frac{d}{dt}V(\gamma_0(s), \frac{d\gamma_0(s) }{ds})<0$ for almost all $s\in[0,1]$ then $\frac{d}{dt}D_\gamma(x_a, x_b)<0$. Then a let $d(x_a, x_b) = \inf_\gamma D_\gamma(x_a, x_b)$ where the infimum is taken over continuously differentiable paths connecting $x_a$ and $x_b$, as with a Riemannian metric \cite{boothby1986introduction}.

Now consider solutions $\gamma: [0, 1]\times \mathbb R^+\rightarrow K$ of $\dot x= E(x)^{-1}f(x,u)$ with initial conditions $\gamma(s, 0) = \gamma_0(s)$ for each $s\in[0,1]$. Suppose at each point $\gamma(s,t)$ along the path, the system is ``sped up'' or ``slowed down'' by a factor $\dot\tau(s,t)$, which is a differentiable function of $s$. 

Under such time reparametrization, it follows that
\begin{align}
&\frac{d}{dt}V\left(\gamma(s,\tau),\frac{\partial \gamma(s,\tau)}{\partial s}\right)  = \dot\tau\left(\left. \frac{d}{dt}V\left(\gamma(s, \tau),\frac{\partial \gamma(s,\tau)}{\partial s}\right)\right|_{\frac{\partial\dot\tau}{\partial s}=0} \right. \notag \\& \left.+ \frac{\partial\dot\tau}{\partial s}\left.\frac{\partial V(\gamma(s,\tau), \delta_x)}{\partial \delta_x}\right|_{\delta_x = \frac{\partial \gamma(s,\tau)}{\partial s} }E(\gamma(s,\tau))^{-1}f(\gamma(s,\tau),u) \right), \label{Vdot_with_z}
\end{align}
where $\tau=\tau(s,t)$.

If there exists a $\frac{\partial\dot\tau(s,t)}{\partial s}\in \mathbb R$ such that the quantity in \eqref{Vdot_with_z} is negative for all $s\in[0, T]$, then by choice of $\dot\tau(s,t)$ the distance between two points $d(x_a, x_b)$ can be made to decrease, for some time reparametrizations.

Since \eqref{Vdot_with_z} is affine in $\frac{\partial\dot\tau(s,t)}{\partial s}$, a sufficient condition for $V(\gamma(s,\tau),\frac{\partial \gamma(s,\tau)}{\partial s})$ to be decreasing is the following:
\begin{equation}\label{eqn:transcon}
\frac{d}{dt}V(x, \delta_x)<0 \forall x\in K, \delta_x \in \mathbb R^n :\frac{\partial}{\partial \delta_x} V(x, \delta_x)\dot x = 0.
\end{equation}

Furthermore, the particular choice of $x=\gamma(s)$ and $\delta_x = \frac{\partial\gamma}{\partial s}$ implies that $\frac{d}{dt}V(\gamma(s), \frac{d\gamma }{ds})<0$.

In this paper we propose the following form of
\[
V(x,\delta_x) = |E(x) \Pi_{\dot x}\delta_x|_Q^2+|\pi_{\dot x}\delta_x|^2
\]
we have $\frac{\partial}{\partial \delta_x} V(x, \delta_x) = \delta_x' \Pi_{\dot x}'E(x)'QE(x) \Pi_{\dot x}+\delta_x'\pi_{\dot x}'\pi_{\dot x}$ and 
by construction of $\Pi_{\dot x}$ and $\pi_{\dot x}$ it follows that $\frac{\partial}{\partial \delta_x} V(x, \delta_x)\dot x = 0$ if and only if $\delta_x$ is such that $\pi_{\dot x}\delta_x = 0$ and $\Pi_{\dot x}\delta_x = \delta_x$.

For such $\delta_x$, $V(x, \delta_x) = |E(x) \Pi_{\dot x}\delta_x|_Q^2$ and we have
\[
\frac{d}{dt} V(x, \delta_x) = 2\delta_x{\Pi_{\dot x}}' E'Q((F+E\dot \Pi_{\dot x}) \Pi_{\dot x}\delta_x.
\]
Negativity of the right hand side is implied by the conditions of the theorem.

This verifies that any two points $x_a, x_b$ in $K$, there exists monotonic smooth functions $\tau_a(t):[0,\infty)\rightarrow [0, \infty)$ and  $\tau_b(t):[0,\infty)\rightarrow [0, \infty)$ such that $x(\tau_a(t))\rightarrow x(\tau_b(t))$ as $t\rightarrow \infty$, a form of incremental stability sometimes referred to as Zhukovsky stability. I follows from the strict forward-invariance of $K$ that all solutions have an $\Omega$-limit set in $K$, and it follows from the Zhukovsky stability that all solutions have the {\em same} $\Omega$-limit set.

Furthermore, it is known that solutions having this property have an $\Omega$-limit set that is a periodic cycle \cite{yang2000liapunov}, which we denote $x^\star(\cdot)$. This completes the proof of the theorem.
\end{prf}

The conditions in Theorem \ref{thm:contraction} are only approximately imposed by our identification procedure. Let us describe the approximations:
\begin{enumerate}
\item The contraction condition \eqref{eqn:contraction} is defined using $\Pi^r_{\dot x}$ and $\Pi_{\dot x}$ with $\dot x = E(x)^{-1}f(x,u)$ . Finiteness of the TRIE implies a similar condition but defined using $\Pi^r_{\dot{\tilde x}}$ and $\Pi_{\dot{\tilde x}}$, i.e. the projections are defined to be transversal the state derivative $\dot {\tilde x}$ from the {\em data} rather than the {\em model}. The reason for this is that $\Pi^r_{\dot x}$ and $\Pi_{\dot x}$ are highly nonlinear functions of the model parameters, and there does not appear to be a straightforward way to convexify condition \eqref{eqn:contraction} with respect to $e(x)$ and $f(x)$. In contrast,  $\Pi^r_{\dot{\tilde x}}$ and $\Pi_{\dot{\tilde x}}$ can be computed directly from the data in advance of any identification procedure.
\item Theorem \ref{thm:contraction}  assumes that some set $K$ has been proven to be forward invariant, which we do not explicitly impose. It is possible to give convex conditions for boundedness in a similar manner to the global stability constraint in \cite{Tobenkin10a}, however experience has shown that this is rarely necessary.
\item Supposing we did verify existence of a strictly forward-invariant set $K$, the theorem assumes that condition \eqref{eqn:contraction} holds for all $x$ in $K$, whereas we will impose the similar condition only at the points where we have data samples. The reason for this is simply that we cannot construct $\Pi^r_{\dot{\tilde x}}$ and $\Pi_{\dot{\tilde x}}$ elsewhere.
\end{enumerate}

With regards to point 1) above, due to the strict negativity in \eqref{eqn:contraction} and the smoothness all functions, it is clear that the condition will still hold for $\Pi^r_v$ and $\Pi_v$ where $v$ is in some neighborhood of $\dot x$. Thus, if the identification procedure is sufficiently successful and $\dot x \approx \dot {\tilde x}$, then the condition will still hold.

With regards to point 2), it has been observed in experiments that if a sufficiently rich excitation is used so that the data points ``fill'' the state space in the vicinity of a cycle relatively densely (in the non-mathematical sense) then it is very difficult to find a model by the proposed method that does {\em not} achieve a stable limit cycle.

%%% Local Variables: 
%%% mode: latex
%%% TeX-master: "transID_TAC"
%%% End: 

\section{Implementation as a Semidefinite Program}\label{sec:impl}

We now discuss practical considerations for data preparation and
minimization of the upper bounds using semidefinite programming.
\subsection{Extracting States from Input/Output Data}
\label{sec:state}
The RIE formulation assumes access to approximate state observations,
$\tilde x(t)$.  In most cases of interest, the full state of the
system is not directly measurable and extraction of a state vector is a challenging problem in its own right.  In practice, our solutions have
been motivated by the assumption that future output can be
approximated as a function of recent input-output history and future
input.  For autonomous systems this assumption is well motivated by the Takens embedding theorem \cite{takens1981detecting}. A common method for choosing a set of time-delays is to optimize mutual information \cite{Fraser86}. An alternative to pure time delays that we have had success with is a
linear filter bank applied to the output, such as a series of Laguerre filters \cite{wahlberg1991system, Chou99}). This has the advantage that the derivatives of these variables can to be
calculated analytically.

Projection-based methods such as subspace identification \cite{Overschee94} and proper-orthogonal-decomposition (POD) \cite{Holmes98} are common methods for generating a state space. Modifications of POD based on balanced truncation have also been proposed for nonlinear model reduction \cite{lall2002subspace} However, even in fairly benign cases of nonlinear systems, one expects the input-output histories to
live near a nonlinear submanifold of the space of possible histories. Connections
between nonlinear dimensionality reduction and system identification
have been explored in some recent papers, e.g.
\cite{Rahimi07} and \cite{Ohlsson08}.

The above procedures involve first estimating a set of states, and then performing the identification. However, there are alternatives to these two-step procedures. In \cite{ramsay2007parameter}, a non-parametric smoothing spline was combined in a nested minimization with an equation-error-based parameter identification. In \cite{schon2011system} the Expectation-Maximization procedure was applied, resulting in a successive alternation between state estimation via a particle filter and parameter estimation via maximum likelihood.

Each of the above procedures consists of some combination a method of deriving states is combined in some way with a method of approximating dynamics: usually either by nonlinear programming on simulation error (or maximum likelihood) or basic least squares on the equation error. In all these cases, the TRIE approach presents an alternative procedure for the approximation of the dynamics.

\subsection{Semidefinite/Sum-of-Squares Programming Formulation}
For each data point, the upper bound on the local TRIE
is the supremum of a concave quadratic form in $\Delta$.  So long as
$e,f$ and $g$ are chosen to be linear in the decision variables, this upper
bound can be minimized by introducing a linear matrix inequality (LMI) for each of $N$ data-points:  we introduce a slack variable $s_i$ for each
data-point $ s_i \ge \hat{\cl E}^\perp_Q(t_i)$ and minimize their sum.

We parametrize $e, f,$ and $g$ as polynomials, so that a sum-of-squares relaxation \cite{Parrilo03a} is used to enforce the well-posedness constraint on $E(x)$, resulting in the following convex optimization problem, 
\begin{align}
&\sum_i s_i \rightarrow \min \textrm{ subject to} \notag\\
&\begin{bmatrix}s_i-\frac{1}{2}|\Pi^r\Delta|_{Q^{-1}}^2+(\Pi^r\Delta)'\Delta_e^- \notag
  - |\Delta_y|^2 & {\Delta_e^+}'\\{\Delta_e^+}&2Q^{-1}\end{bmatrix}>0, \\& \forall \Delta\in\mathbb R^{n-1}, i=1, 2, ..., N, \notag\\
&E(x)+E(x)'\ge I \ \forall x\in\mathbb R^n, \notag
\end{align}
which is easily implemented as a semidefinite program using software such as SPOT \cite{megretski2010spot} or Yalmip \cite{lofberg2004yalmip}.

%%% Local Variables: 
%%% mode: latex
%%% TeX-master: "transID_TAC"
%%% End: 

\section{Experimental Results on Live Neurons}\label{sec:exp}

We now demonstrate the method by identifying the membrane-potential
dynamics of a live, {\it in-vitro} hippocampal neuron. A micropipette
is used to establish an interface with the cell such that current can
be injected into the soma, and the membrane potential response can be
recorded. A microscopic photograph of the neuron and patch is shown in
Figure~\ref{fig:neuron}. The preparation of the culture is described
in the appendix.  

The membrane dynamics of a single neuron are highly complex: at low currents the system responds like a low-order passive linear system. However, when certain ion channels are activated rapid {\em spiking} can occur. After spiking there is a {\em refractory} period in which sensitivity is reduced, and with sufficient current input spiking can repeat at an input-dependent frequency.

\begin{figure} \centering \includegraphics[width=\columnwidth]{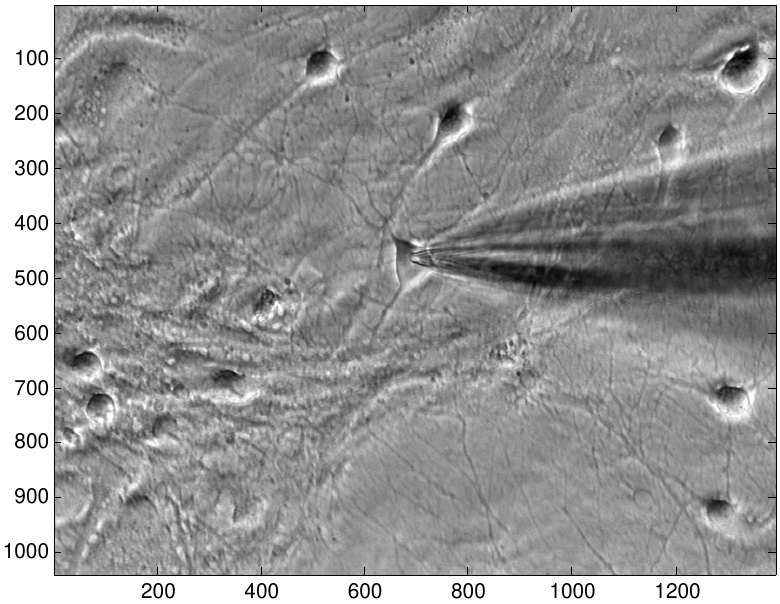}
  \caption{The neuron in culture and the glass micropipette electrode used to interface to it. Imaging: phase contrast image at 20$\times$ magnification on an inverted Olympus IX-71 microscope. Scale: 100 pixels (marked on axes) = 43 microns.}\label{fig:neuron} \end{figure}

There is a spectrum of models of neuron dynamics, ranging from simple ``integrate and fire'' models to highly complex biophysical models of ion channels and conductances \cite{Herz06}.  Threshold based models
generally have a very small number of parameters, but do not provide
high fidelity reproduction of the membrane potential dynamics.  By
contrast biophysical models can be very accurate, but are highly nonlinear and
are very difficult to identify \cite{Geit08} -- they  can have many locally optimal fits in disconnected regions of parameter space
\cite{Achard06}.  In this section we use the proposed method to identify a black-box
nonlinear model with comparatively few states (three) which reproduces
the experimentally observed spiking and subthreshold behavior with very
high fidelity.

\subsection{Identification Results}

The models we consider were of the form \eqref{eqn:model_ef}, \eqref{eqn:model_g} with three states. Each element of the matrix $e(x)$ was a third-degree polynomial in $x$, and each element of $f(x,u)$ was third degree in $x$ and affine in $u$, and $g(x,u)$ was affine in $x$ and $u$. The number of decision variables were 60 for $e(x)$, 240 for $f(x,u)$, 7 for $g(x,u)$ and 6 for $Q$, giving a total of 313 decision variables to define the model and storage function.

As discussed in Section~\ref{sec:impl}, we must find a good proxy
for the internal state of the system.  Here we used two Laguerre
filters with identical pole locations to summarize the recent history
voltage history.  All signals were smoothed with a simple nonparametric smoother.

In the first set of results we show, three increasing step currents are applied to
the neuron resulting in increasing firing rate and a characteristic
change in the spike amplitude and shape.

\begin{figure} \centering \includegraphics[width=\columnwidth]{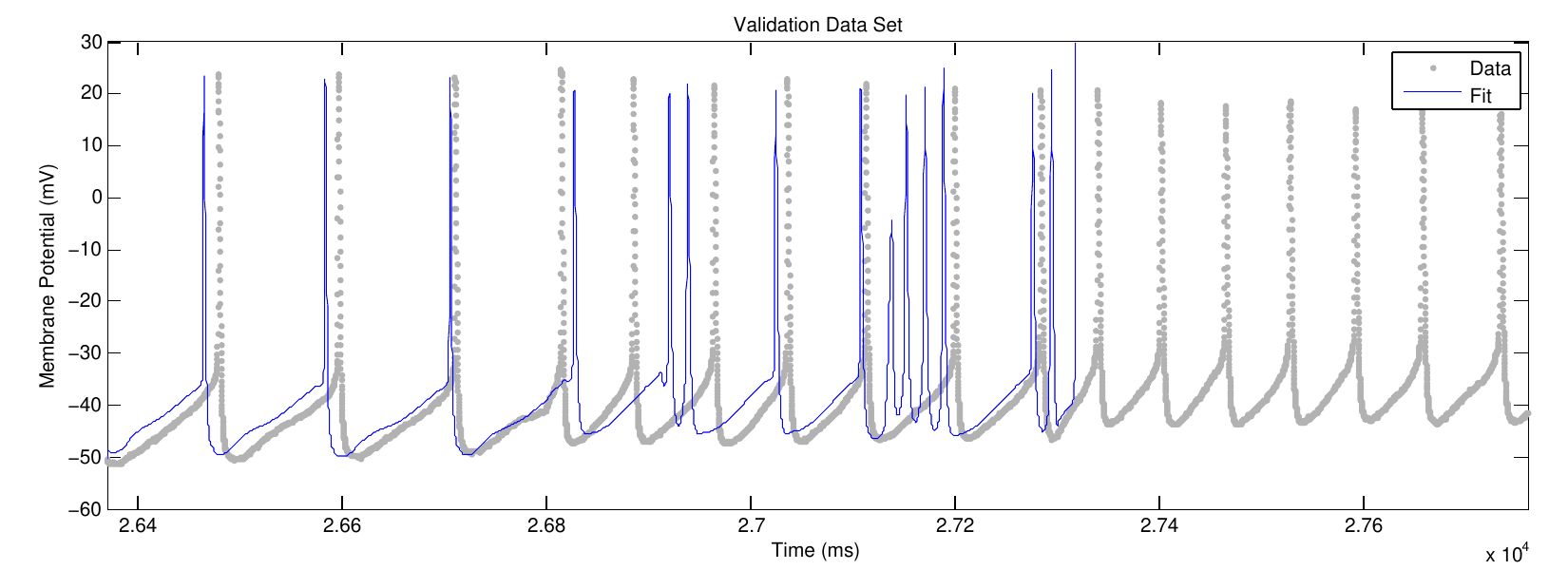}
\includegraphics[width=\columnwidth]{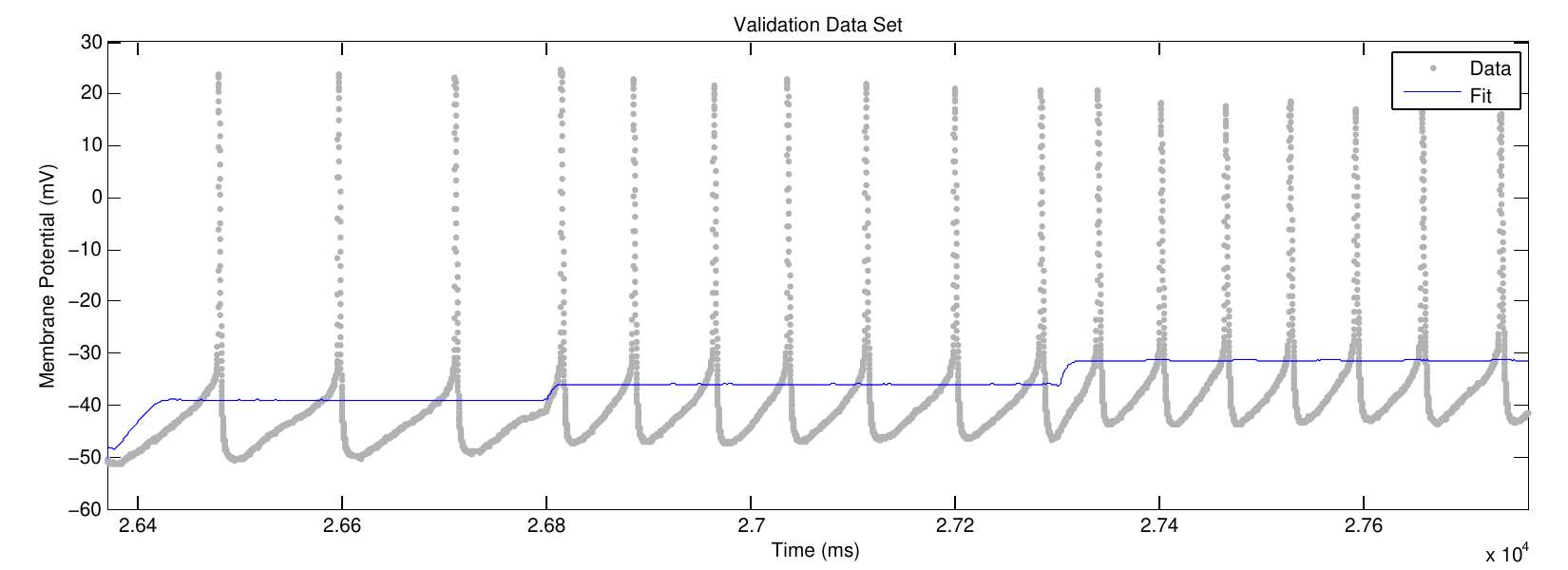}
\includegraphics[width=\columnwidth]{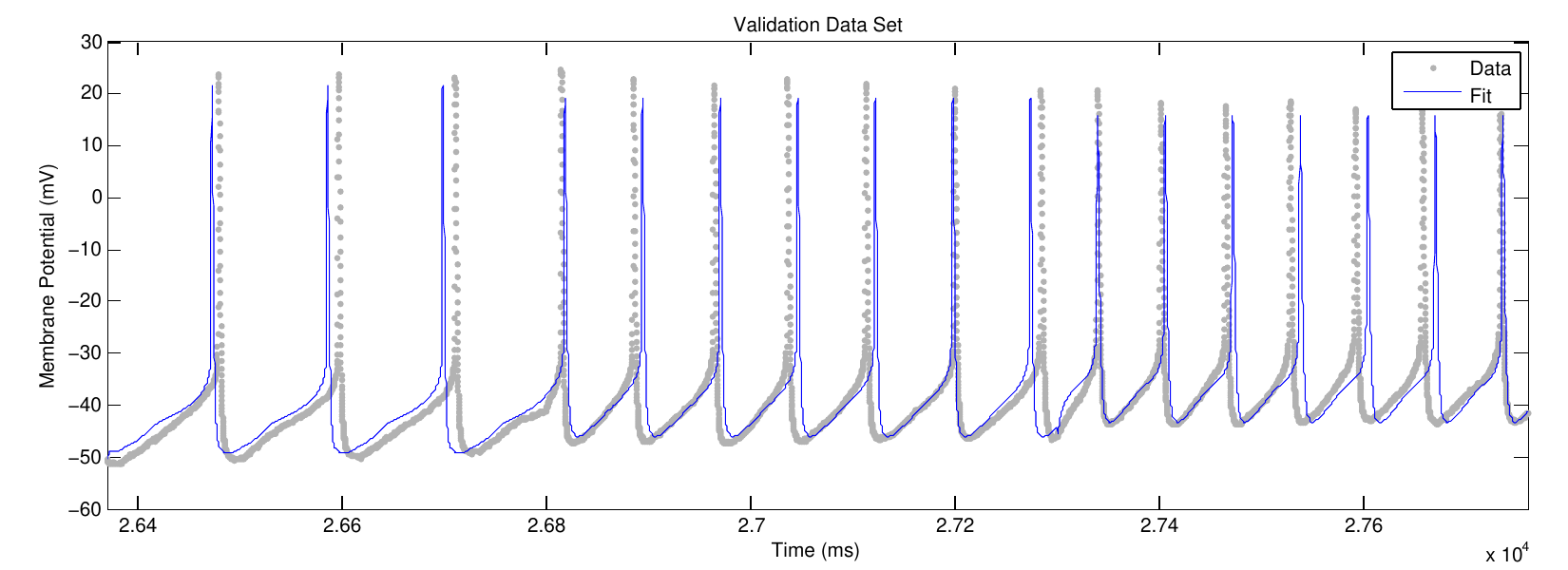} \caption{A neuron was subjected to increasing
step-currents, and spiked with increasing frequency.  Long-term simulation of an equation error fit (top) is unstable.  RIE (middle) minimization provides an overly stable fit.  The
proposed TRIE method (bottom) reproduces the spikes accurately.}\label{fig:neuron-fit} \end{figure}

Figure~\ref{fig:neuron-fit} presents a comparison of fit performance
using three methods. The first is equation error minimization, i.e. simply optimizing
\[
\sum_i|E(\tilde x (t_i))\dot {\tilde x}(t_i) - f(\tilde x (t_i),\tilde u (t_i))|^2\rightarrow \min
\]
subject to the well-posedness constraint $E(x)+E(x)'\ge I$ but without
constraints on stability or long-term simulation error (this is similar in principle to NARX and prediction error methods). The second method is the comparison is the original RIE minimization from \cite{Tobenkin10a}, and the third is
the proposed Transverse RIE method.

We see that while equation error minimization (top)
leads to initially good performance, the model simulation diverges sharply from the recording at about 2.73ms.
Fitting with the RIE (middle) leads to the anticipated overly stable
model dynamics (see Section \ref{sec:rie}), for each level of input the model converges to an equilibrium ``averaging'' the output levels.  The final plot presents
the Transverse RIE identification, which matches the experimentally observed spike patterns very well.

Note that, although the main theoretical results only apply for constant inputs, here we see that the method also works well for transients between piecewise-constant inputs.

Another point to note is that, for the TRIE fit, most of the spikes occur at {\em nearly but not exactly} the right time. It is clear that if one considered simulation error $\int_0^T|y(t)-\tilde y(t)|^2dt$ without adjusting for phase, there would be substantial errors recorded in the intervals where either the model or the true neuron has spiked without the other. This further motivates the concept of ``orbital simulation error''.

\begin{figure} \centering \includegraphics[width=\columnwidth]{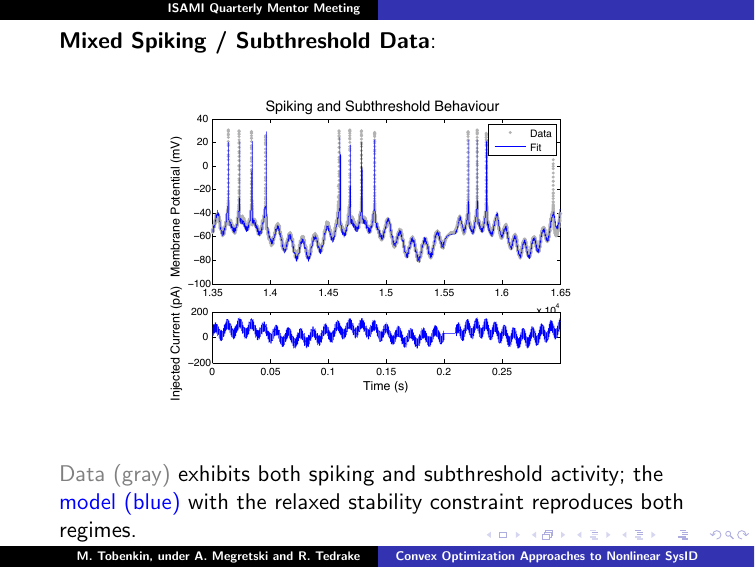}
\caption{Response of real neuron and TRIE model simulation, showing both stable sub-threshold behavior and spiking.}\label{fig:neuron-fit2} \end{figure}

We have also had success identifying behavior which
covers both the subthreshold and spiking regime of a neuron with more complex inputs. The
applied stimulus was a variety of multisine signals.
Figure~\ref{fig:neuron-fit2} presents validation of a Transverse RIE
fit on held-out data.  The lower plot is the multisine input in
pico-Amperes.  The upper plot presents the original data and fit.
Both the subthreshold regime and spikes are generally well reproduced. This illustrates that good fits can be achieved with very complex inputs.

%%% Local Variables: 
%%% mode: latex
%%% TeX-master: "transID_cdc"
%%% End: 

\section{Conclusions}\label{sec:conc}

This paper has introduced a new technique -- Transverse Robust
Identification Error -- for identification of
nonlinear systems which may produce autonomous oscillations,
i.e. system oscillations which are produced internally by the dynamics
rather than as a response to a periodic input.

A convex optimization procedure is developed which minimizes an upper
bound on a local measure of simulation error -- the long-term
divergence between a model simulation and the recorded data.

A theorem was proved giving conditions for a model to have a unique
stable limit cycle, and it was shown to be closely related to the
conditions that are imposed by TRIE.

The proposed method worked well on the challenging problem
of accurately modeling the membrane dynamics of a live neuron from experimental current-voltage
recordings. The input-dependent absence or presence, and frequency of
repetition, of spiking events was well captured in the model.

The method can be implemented on general purpose semidefinite
programming solvers. However, to do so means introducing a large
number of slack variables and LMI constraints. Future work will include investigating dedicated solvers and methods for
extracting states, as well as testing the proposed
method on a wider range of applications.

%%% Local Variables: 
%%% mode: latex
%%% TeX-master: "transID_cdc"
%%% End: 

\appendix

\subsection{Live Neuron Experimental Procedure}

Primary rat hippocampal cultures were prepared from P1 rat pups, in
accordance with the MIT Committee on Animal Care policies for the
humane treatment of animals. Dissection and dissociation of rat hippocampi
were performed in a similar fashion to \cite{Hagler01}. Dissociated
neurons were plated at a density of 200K cells/mL on 12 mm round glass
coverslips coated with 0.5 mg/mL rat tail collagen I (BD Biosciences)
and 4 $\mu$g/mL poly-D-lysine (Sigma) in 24-well plates. After 2
days, 20 $\mu$M Ara-C (Sigma) was added to prevent further growth
of glia. 

Cultures were used for patch clamp recording after 10 days in vitro.
Patch recording solutions were previously described in \cite{Bi98a}.
Glass pipette electrode resistance ranged from 2-4 M$\Omega$.
Recordings were established by forming a G$\Omega$ seal between
the tip of the pipette and the neuron membrane. Perforation of the
neuron membrane by amphotericin-B (300 $\mu$g/mL) typically occurred
within 5 minutes, with resulting access resistance in the range of
10-20 M$\Omega$. Recordings with leak currents smaller than -100
pA were selected for analysis. Leak current was measured as the current
required to voltage clamp the neuron at -70 mV. Synaptic activity
was blocked with the addition of 10 $\mu$M CNQX, 100 $\mu$M APV,
and 10 $\mu$M bicuculline to the bath saline. 
%
% Signals to and from the patch clamp amplifier (Axon Instruments Multiclamp
% 700A) in current clamp mode were sampled at 10 kHz using a National
% Instruments data acquisition card (PCI-6052E).
Holding current was
applied as necessary to compensate for leak current.

\subsection{A Technical Lemma}
% \begin{lem} \label{lem:barDelta} Given the definitions of $\Delta$ in \eqref{eqn:Delta_t}, $\bar\Delta$ in \eqref{eqn:Deltabar_t}, and $\Pi(t)$ in \eqref{eqn:pi}, thenfor all $t\in[0,T]$:
% $\bar\Delta(t) = \Pi(t)\Delta(t)$
% \end{lem}
% \begin{prf}
% By definition,
% \[
% \bar\Delta(t) := \lim_{\theta\rightarrow 0}\frac{1}{\theta}(x_\theta(t) - \tilde x(\tau_\theta))
% \]
% where
% \[
% \tau_\theta = \arg\min_{\tau\in [0, T]}|\tilde x(\tau)-x_\theta|.
% \]
% First note that for any $\theta$ sufficiently small,
% \[
% \Pi(\tau)\frac{x_\theta(t) - \tilde x(\tau)}{\theta} = \frac{x_\theta(t) - \tilde x(\tau)}{\theta}
% \]
% since $\tilde x(\tau)$ is the closest point to the solution set $\{\tilde x(\cdot)\}$ to $x_\theta(t)$, and so $x_\theta(t)-\tilde x(\tau)$ is in the subspace orthogonal to $\dot {\tilde x}(\tau)$ and is therefore invariant under the projection $\Pi(\tau)$. And since $\Pi(\tau)\rightarrow\Pi(t)$ as $\theta\rightarrow 0$ we have
% \[
% \bar\Delta(t) = \Pi(t)\lim_{\theta\rightarrow 0}\frac{1}{\theta}(x_\theta(t) - \tilde x(\tau_\theta)).
% \]

% Now,  define $\bar x_\theta(t) = \tilde x(\tau_\theta)-\tilde x(t)$. Since $\tau_\theta- t= O(\theta)$, it follows that $\bar x_\theta(t)/\theta$ is a rescaling of $\dot{ \tilde x}$ for $\theta$ small, i.e.
% \[
% \lim_{\theta \rightarrow 0}\frac{\bar x_\theta(t)}{\theta} = \alpha \dot{ \tilde x}(t)
% \]
% for some $\alpha\in\mathbb R$.
% Now, by definition $\Pi(t) \dot{ \tilde x}(t) = 0$ so 
% \[
% \bar\Delta(t) = \Pi(t)\lim_{\theta\rightarrow 0}\frac{1}{\theta}(x_\theta(t) - \tilde x(t)) = \Pi(t)\Delta(t)
% \]
% \end{prf}

\begin{lem}\label{lem:Vdot}
Given the storage function \eqref{eqn:trans_storage}, with $\Delta(t,\tau)\in S(t)$ and $\dot\tau = 1$, then
\begin{equation}\label{eqn:vDbar_dot}
\frac{d}{dt}V(\Delta(t),t)=2\Delta{\Pi}' E'Q((F+E\dot \Pi) \Pi\Delta+\epsilon_x).
\end{equation}
\end{lem}
\begin{prf}
% \[
% \dot {V} = 2\Delta{\Pi}' E'Q((F+E\dot \Pi) \Pi\Delta+\epsilon_x)
% \]
Let
\[
v = E(\tilde x(t)) \Pi(t)\Delta, \delta = \pi(t)\Delta,
\]
so  $V(\Delta, t) = v'Qv+\delta'\delta$ and 
\[
\dot {V} = 2v'Q\dot v+2\delta\dot\delta
\]
but we have $\Pi(t)\Delta(t)=\Delta(t)$ and $\pi(t)\Delta(t) = 0$, so 
\[
\dot {V} = v'Q\dot v.
\]
We now derive an expression for $\dot v$.

Let $t_1=t$ and decompose $\Delta(t)$ into two components $\Delta_\perp(t)+\bar\delta(t)$,
transversal and tangential, with respect to $\dot x(t_1)$:
\[
\Delta_\perp(t) = \Pi(t_1)\Delta(t), \bar\delta(t)= \pi(t_1)\Delta(t)
\]
Note that
this decomposition is based on the transversal and tangential decomposition at a {\em fixed}
time $t_1$, not based on a rotating coordinate system. By the chain rule,
\begin{eqnarray}
\dot v&:=&\left.\frac{d}{dt}\left[E(x(t))\Pi(t)\Delta(t))\right]\right|_{t=t_1}\notag\\
&=&\left.\frac{d}{dt}\left[E(x(t))\Pi(t_1)\Delta(t)+E(x(t_1))\Pi(t)\Delta(t_1)\right]
\right|_{t=t_1}\notag\\
&=&\left.\frac{d}{dt}[E(x(t))\Pi(t_1)(\Delta_\perp(t)+\bar\delta(t))]\right|_{t=t_1}\notag\\
&&+E(x(t_1))\dot\Pi(t_1)\Delta(t_1)\notag
\end{eqnarray}
but by definition $\Pi(t_1)\bar\delta(t)=0$ and $\Pi(t_1)\Delta_\perp(t) = \Delta_\perp(t)$ for all $t$, so
\begin{eqnarray}
\dot v&=&\left.\frac{d}{dt}[E(x(t))\Delta_\perp(t)]\right|_{t=t_1}+E(x(t_1))\dot\Pi(t_1)\Delta(t_1)+\notag\\
&=&F(x(t_1),u(t_1))\Delta_\perp(t_1) +\epsilon_x+E(x(t_1))\dot\Pi(t_1)\Delta(t_1)\notag
\end{eqnarray}
But $\Delta_\perp(t_1) = \Delta(t_1)=\Pi(t_1)\Delta(t_1)$, so for this particular $\Delta$, 
\begin{equation}\label{eqn:Vdoteq}
\dot {V}(\Delta(t), t) = 2\Delta{\Pi}' E'Q((F+E\dot \Pi) \Pi\Delta+\epsilon_x)
\end{equation}
which proves the lemma.
\end{prf}

%%% Local Variables: 
%%% mode: latex
%%% TeX-master: "transID_cdc"
%%% End: 

\bibliographystyle{IEEEtran}
\bibliography{elib}

\end{document}